\documentclass[english]{article}
\usepackage[T1]{fontenc}
\usepackage[latin9]{inputenc}
\usepackage{mathrsfs}
\usepackage{amstext}
\usepackage{amssymb}
\usepackage{babel}
\begin{document}
\title{Finite regular semigroups with permutations that map elements to inverses}
\author{Peter M. Higgins }
\date{{}}
\maketitle
\begin{abstract}
We give an account on what is known on the subject of \emph{permutation
matchings}, which are bijections of a finite regular semigroup that
map each element to one of its inverses. This includes partial solutions
to some open questions, including a related novel combinatorial problem.
\end{abstract}

\section{Introduction and background}

The author published three papers {[}3{]}, {[}4{]}, and {[}5{]} in
Semigroup Forum relating to the question of when a finite regular
semigroup $S$ has a bijection $\phi$ such that $\phi(a)\in V(a)$,
which is to say that $\phi$ maps elements to their inverses. Some
simple observations will natually occur to the reader: if $S$ is
an inverse semigroup, then there is a unique such permutation $\phi$,
and if $S$ is a union of groups, the mapping $a\mapsto a^{-1}$,
where $a^{-1}$ denotes the group inverse of $a$ in $S$, defines
a \emph{permutation matching}, (we shall sometimes abbreviate the
term to \emph{matching}) as these bijections are called. Rectangular
bands are exactly the semigroups where every bijection of $S$ is
a matching, as they are characterised by the identity $x=xyx$, while
the identity mapping is a matching if and only if $S$ satisfies the
identity $x=x^{3}$. 

It does not take long however to locate examples of finite regular
semigroups that have no matching. A minimal example is given by the
7-element orthodox 0-rectangular band $B=\{(i,j):1\leq i\leq2,1\leq j\leq3\}\cup\{0\}$,
where $E(B)=\{(1,2),(1,3),(2,1)\}\cup\{0\}$. The pairs $(2,2)$ and
$(2,3)$ have only one inverse element between them, that being $(1,1)$,
and so no permutation matching exists for $B$. Proposition 1.5 of
{[}3{]} gives an example of an orthodox $0$\emph{-}rectangular band,
with an involution matching, but which has $B$ as a retract. It follows
from this that the class of semigroups that possess a permutation
matching is not closed under the taking of regular subsemigroups,
nor the taking of homomorphic images. The class is closed however
under the taking of direct products.

The nature of the example represented by $B$ suggests that the existence
of permutation matchings will correspond to the satisfaction of Hall's
condition in the bipartite graph featuring two copies of $S$ with
edges joining mutually inverse pairs. Indeed we see this as part (iii)
of the following theorem. 

\vspace{.3cm}

\textbf{Theorem 1.1} {[}3, Theorem 1.6{]} For a finite regular semigroup
$S$ the following are equivalent: 

(i) $S$ has a permutation matching; 

(ii) $S$ is a transversal of $\{V(a)\}_{a\in S}$; 

(iii) $|A|\leq|V(A)|$ for all $A\subseteq S$; 

(iv) $S$ has a permutation matching $(\cdot')$ that preserves the
$\mathscr{H}$-relation; (meaning that $a\,\mathscr{H}\,b\rightarrow a'\,\mathscr{H}\,b'$); 

(v) each principal factor $D_{a}\cup\{0\}$ $(a\in S)$ has a permutation
matching; 

(vi) each $0$-rectangular band $(D_{a}\cup\{0\})/{\cal \mathscr{H}}$
has a permutation matching. 

\vspace{.3cm}

This theorem was applied to show that the full transformation semigroup
$T_{n}$ has a matching {[}3, Theorem 2.12{]}. The argument also extends
to the partial transformation semigroup, $PT_{n}$. 

The proof shows that each $\mathscr{D}$-class $D$ of $T_{n}$ has
a matching, and this was done by partitioning $D$ into so-called
${\cal Q}$ classes, which are themselves unions of $\mathscr{R}$-classes.
To each $\mathscr{R}$-class $R\subseteq D$, (remembering that $\alpha\,\mathscr{R\,\beta}$
in $T_{n}$ if and only if ker$\alpha=$ ker$\beta$), we associate
an increasing sequence of positive integers $s_{R}=(k_{1},k_{2},\cdots,k_{t})$
where the $k_{i}$ run over the cardinalities of the kernel classes
of the members of $R$, so that $t$ is the rank of the members of
$D$. The ${\cal Q}$-class defined by $R$ is then the union of all
$\mathscr{R}$-classes $R'$ such that $s_{R'}=s_{R}$. 

We then show the existence of a permutation matching for each ${\cal Q}$-class
of $T_{n}$. To do this, we consider the bipartite graph $G$ consisting
of two copies, ${\cal Q}$ and ${\cal Q}'$, of the members of ${\cal Q}$,
where $ab'$ is an edge in $G$ if (and only if) $a$ and $b$ are
mutual inverses in $T_{n}$. It is then proved that $G$ is a regular
graph, which is to say each vertex is of a common degree $m\geq1$,
whence each member of ${\cal Q}$ has exactly $m$ inverses in ${\cal Q}$.
It then follows by a standard result from graph theory that $G$ has
a \emph{perfect matching}, which is a set of disjoint edges covering
all vertices of $G$. 

This perfect matching may then be used to generate a permutation matching
of ${\cal Q}$ as follows. Begin with any vertex $a_{1}$ of ${\cal Q}$
and consider the sequence of vertices $a_{1},b_{1}',b_{1}:=a_{2},b_{2}',b_{2}:=a_{3},\cdots$
, which will eventually yields a cycle $a_{1},a_{2},\cdots,a_{r}$
say, where $a_{i}$ and $a_{i+1}$ are mutual inverses (and where
we take $r+1=1)$. Since the edges of $G$ form a perfect matching,
no vertex $a_{i}$ of ${\cal Q}$ can be followed in the ${\cal Q}$-sequence
by a previous vertex of ${\cal Q}$, other than $a_{1}$, as that
would result in the contradiction of two edges of the perfect matching
having a common vertex in ${\cal Q}'$. We then repeat the process
for a new vertex not contained in any of the previous cycles of inverses,
culminating in a permutation matching of ${\cal Q}$.

However the following question remains open.

\vspace{.3cm}

\textbf{Question 1. }Does the semigroup $O_{n}$ of all order-preserving
mappings on a finite $n$-chain have a permutation matching?

\section{Involution matchings}

The unique matching of an inverse semigroup is not just a permutation
but an involution, and likewise may be said for completely regular
semigroups through the matching of elements to their group inverses.
Another example type with a natural \emph{involution matching }as
we shall call them is the class ${\cal OP}_{n}$ of all semigroups
of orientation-preserving mappings (Example 2.6 of {[}3{]}), and the
technique extends to the class ${\cal P}_{n}$ of all semigroups of
orientation-preserving and orientation-reversing mappings (Section
3.2 of {[}4{]}). For background on these classes of semigroups see
{[}2{]}, {[}6{]}, and {[}7{]}. A natural question is then:

\vspace{.3cm} 

\textbf{Question 2. }Does $T_{n}$ have an involution matching? 

\vspace{.3cm}

By a result of Schein {[}8{]}, $T_{n}$ is covered by its inverse
subsemigroups. If $T_{n}$ possessed an inverse cover with the additional
property that the non-empty intersection of any pair of subsemigroups
in the cover were regular (and hence inverse), it would be possible
to use that cover to construct an involution matching. However, it
was shown in {[}5, Theorem 4.2.3 {]} that for $n\geq4$, no such cover
exists. Moreover for $n\geq8$, $T_{n}$ does not have a permutation
matching in which each element is mapped to a \emph{strong inverse
$a'$, }by which we mean an inverse with the property that $\langle a,a'\rangle$
is an inverse subsemigroup of $T_{n}$. 

The case of $T_{n}$ is however only a particular instance of the
general question, which has not been resolved.

\vspace{.3cm}

\textbf{Question 3. }Does every finite regular semigroup with a permutation
matching possess an involution matching?

\vspace{.3cm}

In {[}4, Theorems 3.7 and 3.1{]} it was proved that any (finite) orthodox
semigroup with a permutation matching has an involution matching.
Moreover an orthodox semigroup possesses a permutation matching if
and only if each $0$-rectangular band $B=(D_{a}\cup\{0\})/\mathscr{H}$
$(a\in S)$ has the property that the maximal rectangular subbands
of $B$ are pairwise \emph{similar}, meaning that the ratio of their
number of $\mathscr{L}$- to $\mathscr{R}$ -classes is a constant
throughout $D_{a}$. A wider class of regular semigroups is defined
by the condition on the arrangement of idempotents within ${\cal \mathscr{D}}$-classes:
for all idempotents $e,f,g\in E(S)$, 
\[
e{\cal \,\mathscr{L}}\,f\,\mathscr{R}g\Rightarrow\exists h\in E(S)\,\text{such that\,\ensuremath{e{\cal \,\mathscr{R}}\,h\,{\cal \mathscr{L}\,}g}},
\]
giving an egg-box picture where the four idempotents in question form
a `solid' rectangle within their ${\cal \mathscr{D}}$-class, which
lead to the name \emph{E-solid }semigroups\emph{. }Equivalently, a
regular semigroup $S$ is E-solid if and only if the idempotent generated
subsemigroup $\langle E(S)\rangle$ of $S$ is a union of groups.
Using the fact that each 0-rectangular band $B=(D_{a}\cup\{0\})/\mathscr{H}$
$(a\in S)$ of a finite E-solid semigroup is orthodox, we may extend
the previous results from orthodox to E-solid semigroups {[}5, Theorem
1.3.5{]}.

We can attempt to build an involution from a given permutation matching
as follows. Consider the graph $G(S)$ with vertices corresponding
to the members of $S$ and edges joining each pair of mutual inverses.
Next suppose that $S$ possesses a permutation matching $\phi$ and
consider the subgraph of $G(S)$ consisting of the vertices and edges
of a permutation matching, which is a set of disjoint cycles. Any
even cycle may then be split into mutually inverse pairs. On the other
hand if an odd cycle contains an idempotent, $e$, we may pair $e$
with itself and pair the remaining members of this odd cycle in mutually
inverse pairs. It follows that if $S$ has a permutation matching,
but no involution matching, then \emph{every} permutation matching
of $S$ must contain an odd cycle that is idempotent-free. 

Suppose then that $S$ is an example of a semigroup of minimum cardinality
that has a matching but no involution matching. Then every principal
factor $D_{a}\cup\{0\}$ of $S$ (or simply $D_{a}$ if $D_{a}$ is
the minimum ideal of $S$) has a permutation matching, mapping $D_{a}$
onto itself. If each of these principal factors had an involution
matching, then the union of these involution matchings on the ${\cal \mathscr{D}}$-classes
$D_{a}$ would define an involution matching of $S$. Since we are
assuming that no such involution matching exists, it follows that
some principal factor has a permutation matching but no involution
matching. Therefore by minimality of the cardinal of $S$, it follows
that $S$ is a completely $0$-simple semigroup. (Every completely
simple semigroup, being a union of groups, has an involution matching.)

By Theorem 1.1(vi), we have that $B=S/\mathscr{H}$ also has a permutation
matching. If $S$ contained proper subgroups, it would follow by the
minimality of the cardinal of $S$ that $B$ possessed an involution
matching $\phi$. However this would allow the construction of an
involution matching $\overline{\phi}$ on $S$ as follows. 

We identify the non-zero members of $B$ with the $\mathscr{H}$-classes
$H_{a}$ of $S$. We then define $a\overline{\phi}=a'$, where $a'$
is the unique inverse of $a$ in $H_{a}\phi$. Then $a'\overline{\phi}$
is the unique inverse of $a'$ in $H_{a'}\phi$, and since $H_{a'}\phi=H_{a}$,
we infer that $\overline{\phi}$ would indeed be an involution matching
on $S$. Since this is a contradiction, we reach the conclusion that
any minimal example $S$ is a $0$-rectangular band, thereby establishing
the following fact. 

\vspace{.3cm}

\textbf{Proposition 2.1 }Suppose that $S$ is a finite regular semigroup
of minimum cardinality with the property that $S$ possesses a permutation
matching but no involution matching. Then $S$ is a $0$-rectangular
band.

\vspace{.3cm}

We thus are led to consider $0$-rectangular bands $S$. Let the non-zero
$\mathscr{D}$-class of $S$ have $m$ rows and $n$ columns, where
without loss we may take $m\leq n$. Indeed we also have $m\geq2$,
as if $m=1$ we have a right zero semigroup, which does have involution
matchings. We are led to the question:

\vspace{.3cm}

\textbf{Question 4. }Does there exist a $0$-rectangular band $S$
that has a permutation matching but no involution matching?

\vspace{.3cm}

We may prove that if this is the case then $m$ cannot be a divsor
of $n$, and so in particular our major $\mathscr{D}$-class $D$
of $S$ cannot be square. In what follows we let $R$ and $C$ respectively
stand for the set of $m$ rows and $n$ columns of $D$. Let us write
$n=am$. 

\newpage

\textbf{Lemma 2.2 }Suppose that $S$ has a permutation matching $\phi$.
Then for any collection $T$ of $t$ rows of $D$ $(1\leq t\leq m)$
\begin{equation}
|\{c\in C:c\cap E(T)\neq\emptyset\}|\geq ta.
\end{equation}
And for any collection $T$ of $t$ columns of $D$ $(1\leq t\leq n)$
\begin{equation}
|\{r\in R:r\cap E(T)\}|\geq\frac{t}{a}.
\end{equation}

\vspace{.2cm}

\textbf{Proof} The $t$ rows contain $nt$ members. These are mapped
by $\phi$ into $s$ columns say, which collectively contain $ms$
members. Hence $nt\leq ms$ and so 
\[
s\geq t\cdot\frac{n}{m}=t\cdot\frac{am}{m}=ta.
\]
Since the set on the left hand side of (1) contains the set of $s$
columns, the inequality (1) now follows. 

Similarly, any set of $t$ columns $(1\leq t\leq n)$ is mapped by
$\phi$ into $s$ rows of $D$ so that the $t$ columns collectively
contain $mt$ members. Hence $mt\leq ns$ and so 
\[
s\geq t\cdot\frac{m}{n}=t\cdot\frac{m}{am}=\frac{t}{a},
\]
and then (2) follows. $\Box$

\vspace{.3cm}

\textbf{Proposition 2.3 }Let $S$ be a $0$-rectangular band with
$m$ rows and $n=am$ $(a\geq1)$ columns in its non-zero $\mathscr{D}$-class.
If $S$ has a permutation matching $\phi$, then $S$ has an involution
matching. 

\vspace{.3cm}

\textbf{Proof} Consider the bipartite graph $G$ with vertex set $R\,\dot{\cup}\,C$,
the disjoint union of the respective sets of rows and columns of $D$,
with an edge connecting a row and a column if that row and column
intersect in an idempotent. Note in particular that if $a\phi=b$,
then an edge connects $R_{a}$ to $L_{b}$ (and $R_{b}$ to $L_{a}$). 

From (1) and (2) we see that any set of $t$ rows of $R$ is collectively
adjacent to at least $at$ columns of $C$ in $D$, and any set of
$t$ columns is collectively adjacent to at least $\frac{t}{a}$ rows
of $D$. By Hall's harem theorem, (see Section III.3 of {[}1{]}) there
exists $a$ injective functions, $\pi_{0},\pi_{1},\cdots,\pi_{a-1}$
from $R$ into $C$ with pairwise disjoint ranges, whose union is
$C$. We label the columns in the range of $\pi_{t}$ as $tm,tm+1,\cdots,(t+1)m-1$. 

We may then define an involution matching on $D$ as follows. Let
$(i,j)\in D$ with $j=tm+r$, where $0\leq r\leq m-1$. We define
\begin{equation}
(i,j)\pi=(j\pi_{t}^{-1},i\pi_{t}).
\end{equation}
Then 
\[
(i,j)\pi^{2}=(j\pi_{t}^{-1},i\pi_{t})\pi=(i\pi_{t}\pi_{t}^{-1},j\pi_{t}^{-1}\pi_{t})=(i,j),
\]
so that $\pi$ defines an involution on $D$ (which extends to $S$
by putting $0\pi=0$). Moreover, by definition of adjacency in $G$
we have that $(i,i\pi_{t})\in E(D)$, and also $(j\pi_{t}^{-1},j)\in E(D)$
as $j\pi_{t}^{-1}\pi_{t}=j$. Hence
\[
(i,j)\,(i,j)\pi\,(i,j)=(i,j)(j\pi_{t}^{-1},i\pi_{t})(i,j)
\]
\[
=(i,i\pi_{t})(i,j)=(i,j),
\]
\[
(j\pi_{t}^{-1},i\pi_{t})(i,j)(j\pi_{t}^{-1},i\pi_{t})=(j\pi_{t}^{-1},j)(j\pi_{t}^{-1},i\pi_{t})=(j\pi_{t}^{-1},i\pi_{t}),
\]
thereby verifying that $(i,j)$ and $(i,j)\pi$ are mutual inverses
in $D$. Therefore $\pi$ is indeed an involution matching of $S$.
$\Box$

\vspace{.3cm}

There have been lively exchanges on a problem related to ours, some
of which have been recorded on Mathoverflow {[}9{]}. The conversation
was initiated by Fedya Petrov who explains that our question, (matching
implies involution matching), would be settled in the affirmative
if the same were true of the following purely combinatorial problem. 

A set of $m$ girls (corresponding to rows) have $mn$ balls so that
each girl has $n$ balls. There are $m$ balls of each of $n$ colours
(the colours corresponding to columns). Two girls may exchange the
balls (1 ball for 1 ball, so each girl still has $n$ balls), but
no ball may participate in more than one exchange. The goal is to
achieve a situation where each girl has balls of all $n$ colours
(and so exactly one of each colour). This question we shall call this
the \emph{Colour alignment problem}.

The participants in the discussion then claim proofs of results that
subsume Proposition 2.3 in that the problem may be solved if $n\equiv0,\pm1,\pm2$
(mod $m$). However, at the time of writing, the general question
remains open. 

\vspace{.3cm}

\textbf{Theorem 2.4 }(F. Petrov, private communication) A positive
solution of the Colour alignment problem implies a positive solution
of the involution matching problem.

\vspace{.3cm}

\textbf{Proof }Suppose that the Colour alignment problem may always
be solved. (The exchanges may include vacuous exchanges where a girl
exchanges a ball with itself; in this way we may claim that in our
solution every ball is exchanged exactly once.) Let $S$ be an $m\times n$
$0$-rectangular band with a permutation matching $\phi$. We shall
denote typical rows and columns of $D$ respectively as $i$ and $j$.
In this way an entry of $D$ may be written as $(i,j)$ or simply
denoted by a single symbol $x,y$ etc. In this latter case, if $x=(i,j)$,
we write $x_{1}=i$ and $x_{2}=j.$ In view of this, in what follows,
notation is less cluttered by writing mappings on the left. 

Let the first co-ordinates be enumerated by the set of $m$ girls
and the second coordinates by the $n$ colours. To girl $a$ we assign
the balls of colour $\phi(a,\alpha){}_{2}$ $(1\leq\alpha\leq n)$.
This assigns $n$ balls to each girl. Moreover, since $\phi$ is bijective,
one ball is assigned for each entry of the array $D$. In particular
there are exactly $m$ balls of each colour. 

We construct our involution matching $\Phi$ on $D$ as follows. Suppose
that there is an exchange between girl $a$ and the ball corresponding
to $\phi(a,\alpha)_{2}$, with girl $b$ and the ball corresponding
to $\phi(b,\beta)_{2}$. We then put
\begin{equation}
\Phi(a,\phi(b,\beta){}_{2})=(b,\phi(a,\alpha)_{2}).
\end{equation}

This is a well-defined function, for if not there would be an equality
of the form $\phi(b,\beta)_{2}=\phi(b',\beta')_{2}$, with girl $a$
and the ball corresponding to $\text{\ensuremath{\phi(a,\alpha')_{2}}}$
say, exchanged with girl $b'$ and ball $\phi(b',\beta')_{2}$. But
then, since $\phi(b,\beta)_{2}=\phi(b',\beta')_{2}$ represents equality
of colours, girl $a$ would finish with two balls of this one colour,
which would contradict our assumption that each girl ends with exactly
one ball of each colour. Hence $\Phi$ is a function. Moreover, by
the definition (4), $\Phi(b,\phi(a,\alpha)_{2})=(a,\phi(b,\beta)_{2})$,
so that $\Phi$ is an involution. Next we check that the domain of
$\Phi$ is the complete array, $D$.

For any girl $a$ and ball $\phi(a,\alpha)_{2}$ held by girl $a$,
there is an exchange with some girl $b$ and ball $\phi(b,\beta)_{2}$
held by girl $b$, yielding the point $(a,\phi(b,\beta)_{2})$ as
a member of dom$(\Phi)$. By the argument of the previous paragraph,
these $n$ points are pairwise distinct, and so represent the complete
row of $D$ defined by $a$. It follows that dom$(\Phi)$ contains
all $mn$ points of $D$, whence $\Phi$ is a well-defined involution
on $D$. 

To conclude the proof, we need only note that the pairs on either
side of (4) are mutual inverses as both $(a,\phi(a,\alpha)_{2})$
and $(b,\phi(b,\beta){}_{2})$ are idempotents. Therefore $\Phi\cup\{(0,0)\}$
is an involution matching for our semigroup $S$. $\Box$

\vspace{.3cm}

School of Mathematics, Statistics and Actuarial Science, University
of Essex, Email address: peteh@essex.ac.uk. ORCID:0000-0002-1198-2578.

\end{document}